\documentclass[11pt,a4paper]{amsart}
\usepackage{amssymb}
\usepackage[arrow,matrix]{xy}

\title{Perverse equivalences, BB-tilting, mutations and applications}
\author{Sefi Ladkani}
\address{Max-Planck-Institut f\"{u}r Mathematik, Vivatsgasse 7, 53111 Bonn,
Germany}
\email{sefil@mpim-bonn.mpg.de}

\thanks{This work was supported by a European Postdoctoral Institute
(EPDI) fellowship.}

\DeclareMathOperator{\add}{add}

\DeclareMathOperator{\End}{End}
\DeclareMathOperator{\Ext}{Ext}
\DeclareMathOperator{\gldim}{gl.\!dim}
\DeclareMathOperator{\h}{H}
\DeclareMathOperator{\Hom}{Hom}
\DeclareMathOperator{\id}{id}
\DeclareMathOperator{\modf}{mod}
\DeclareMathOperator{\pd}{pd}
\DeclareMathOperator{\per}{per}
\DeclareMathOperator{\rad}{rad}
\DeclareMathOperator{\RHom}{\mathbf{R}Hom}
\DeclareMathOperator{\Tr}{Tr}

\newcommand{\BB}{\mathrm{BB}}
\newcommand{\bQ}{\mathbb{Q}}
\newcommand{\bZ}{\mathbb{Z}}
\newcommand{\gL}{\Lambda}

\newcommand{\cA}{\cC_A}
\newcommand{\cB}{\cC_B}
\newcommand{\cC}{\mathcal{C}}
\newcommand{\cD}{\mathcal{D}}
\newcommand{\cE}{\mathcal{E}}
\newcommand{\cS}{\mathcal{S}}

\newcommand{\dA}{\cD^b(A)}
\newcommand{\dB}{\cD^b(B)}

\newcommand{\wt}{\widetilde}

\theoremstyle{plain}
\newtheorem{thm}{Theorem}[section]
\newtheorem{prop}[thm]{Proposition}
\newtheorem{lemma}[thm]{Lemma}
\newtheorem{cor}[thm]{Corollary}

\theoremstyle{definition}
\newtheorem{defn}[thm]{Definition}
\newtheorem{remark}[thm]{Remark}
\newtheorem{example}[thm]{Example}
\newtheorem{quest}[thm]{Question}

\numberwithin{equation}{section}

\begin{document}

\begin{abstract}
We relate the notions of BB-tilting and perverse derived equivalence at
a vertex. Based on these notions, we define mutations of algebras,
leading to derived equivalent ones.

We present applications to endomorphism algebras of cluster-tilting
objects in 2-Calabi-Yau categories and to algebras of global dimension
at most 2.
\end{abstract}

\maketitle

\section*{Introduction}

In~\cite[\S2]{BrennerButler80}, Brenner and Butler introduced a
construction of tilting modules, known as BB-tilting modules,
generalizing APR-tilts~\cite{APR79}, which are themselves
generalizations of the BGP reflection functors introduced
in~\cite{BGP73}. Roughly speaking, for a finite dimensional
algebra $A$ over a field, the associated BB-tilting modules are
parameterized by the simple $A$-modules satisfying certain
homological conditions (to be recalled in
Section~\ref{ssec:BBtilt} below).

On the other hand, Chuang and Rouquier introduced the notion of
perverse Morita equivalences, which are certain derived equivalences
that are controlled by filtrations $\cS = (\phi = \cS_0 \subset \cS_1
\subset \dots \subset \cS_r)$ of the set of isomorphism classes of
simple $A$-modules together with perversity functions $p : \{1, \dots,
r\} \to \bZ$, see~\cite[\S2.6]{Rouquier06}.

In this paper we focus on the special case where the
filtration $\cS$ has only two levels and $|p(2)-p(1)|=1$.
This case leads to complexes of $A$-modules consisting of only two
projective terms concentrated in consecutive degrees. Such complexes
have already been discussed in the literature, see for
example~\cite{HoshinoKato02} and the discussion of the associated
torsion theory in~\cite{HKM02}.

Many of the constructions of tilting complexes introduced in the
literature to show derived equivalences for various kinds of algebras
are in fact such perverse equivalences. Examples are
\cite[Theorem~11.5.1]{Linckelmann98} for symmetric algebras;
\cite[Theorems~29, 32]{AvellaAlaminos07} for gentle algebras;
\cite{Bastian09}, \cite[\S3.1]{BHL09} and~\cite[Theorems~3.12,
3.15]{Murphy08} for cluster-tilted algebras; and \cite{Vitoria09} for
certain Jacobian algebras.

We start in Sections~\ref{sec:pre} and~\ref{sec:tilt} by reviewing the
relations between these two notions of tilting for a finite dimensional
algebra $A$ given by a quiver with relations. To any vertex $k$ without
loops in the quiver, there are two complexes $T^-_k$ and $T^+_k$
attached. We list equivalent formulations, in homological as well as in
combinatorial language, that these complexes are tilting complexes,
hence leading to perverse equivalences with respect to the two-level
filtration $\cS$ with $\cS_1 = \{k\}$. Furthermore, we show that
BB-tilting can always be regarded as a perverse equivalence induced by
a complex $T^-_k$, and conversely, under certain additional
assumptions, the perverse equivalence given by $T^-_k$ is given by
BB-tilting.

This leads us to define, in Section~\ref{sec:algmut}, mutations of
finite dimensional algebras over an algebraically closed field, leading
to derived equivalent ones.
To any
vertex $k$ without loops in the quiver of such an algebra $A$, the
negative mutation $\mu^-_k(A)$ is defined as the endomorphism algebra
of $T^-_k$ when it is a tilting complex. Similarly, the positive
mutation $\mu^+_k(A)$ is defined as the endomorphism algebra of $T^+_k$
when it is a tilting complex. Thus, a vertex $k$ leads to at most two
mutations of $A$, and they are derived equivalent to $A$. In addition,
we define the BB-mutation $\mu^{\BB}_k(A)$ as the endomorphism algebra
of the BB-tilting module associated with $k$, when it is defined.
Hence, when $\mu^{\BB}_k(A)$ is defined, so is $\mu^-_k(A)$ and they
coincide. Analog notions for Calabi-Yau algebras have been introduced
in~\cite{IyamaReiten08}.

From a K-theoretical viewpoint, we also show that the change-of-basis
transformation of the Grothendieck groups induced by a mutation of an
algebra coincides with that occurring in the matrix mutation, in the
sense of Fomin and Zelevinsky~\cite{FominZelevinsky02}, of the
skew-symmetric matrix corresponding to the quiver of the algebra,
provided that it has no loops and 2-cycles.

Mutations of algebras arise naturally when considering endomorphism
algebras of objects related by approximation sequences in additive
categories. This has been essentially observed by Hu and Xi
in~\cite{HuXi08}. In Section~\ref{sec:endo} we show that in an additive
$\Hom$-finite category over an algebraically closed field with split
idempotents, if $\gL$ and $\gL'$ are the endomorphism algebras of two
basic objects related by replacing the indecomposable $k$-th summand by
another through an approximation sequence, then the mere existence of
the two BB-mutations $\mu^{\BB}_k(\gL)$ and $\mu^{\BB}_k(\gL'^{op})$
automatically implies that they take the ``correct'' values, namely
$\gL'$ and $\gL^{op}$, respectively.

Such approximation sequences, known as exchange sequences, appear in
relation with mutations of cluster-tilting objects in 2-Calabi-Yau
(2-CY) categories~\cite{BIRSc09,BMRRT06,GLS06,IyamaYoshino08}, studied
in connection with representation theoretical interpretation of cluster
algebras~\cite{FominZelevinsky02}. In Section~\ref{sec:2CY} we apply
the results of the previous section to the study of endomorphism
algebras of cluster-tilting objects in $\Hom$-finite idempotent split
Frobenius stably 2-CY categories, as well as in such triangulated 2-CY
categories, where these algebras are known as 2-CY-tilted algebras.
The latter algebras include the cluster-tilted
algebras~\cite{BMR07} and more generally~\cite{Amiot08,BIRSm08}
finite dimensional Jacobian algebras of
quivers with potentials~\cite{DWZ08}.

We generalize a result in~\cite{GLS06} and show that in the Frobenius
case, the endomorphism algebra of a cluster-tilting object $U$ admits
all the BB-, negative and positive mutations at any vertex
corresponding to a non projective-injective summand of $U$,
and moreover all these mutations coincide with the endomorphism algebra
of the mutation of $U$ at that summand.

In the triangulated case, the picture is more complicated, and
mutation of cluster-tilting objects does not always lead to a mutation
of their endomorphism algebras.
Indeed, neighboring two 2-CY-tilted algebras, that is, endomorphism
algebras of two cluster-tilting objects related by a mutation, are
always near-Morita equivalent~\cite{BMR07,KellerReiten07,Ringel07},
but not necessarily derived equivalent.
However, several derived equivalence classifications of cluster-tilted
algebras~\cite{Bastian09,BHL09,BuanVatne08} have revealed that there
are far less derived equivalence classes than isomorphism classes of
such algebras (at least when their number is finite). Moreover, these
classifications rely on showing that sufficiently many pairs of
neighboring algebras are in fact derived equivalent.

In an attempt to provide a conceptual explanation of these facts,
we study the conditions that two neighboring 2-CY-tilted algebras $\gL$
and $\gL'$ are related by BB-mutation. Following~\cite{BHL09},
we call in this case the corresponding mutation of the cluster-tilting objects
a good mutation, since it leads to a mutation of their endomorphism algebras.

Obviously, a necessary condition for $\gL' = \mu^{\BB}_k(\gL)$
is that $\mu^{\BB}_k(\gL)$, hence by symmetry also $\mu^{\BB}_k(\gL'^{op})$,
are defined.
We show that this condition is also sufficient, that is,
the existence of these two BB-mutations automatically implies that they
take the ``correct'' values, namely $\gL'$ and $\gL^{op}$,
respectively, yielding a good mutation between $\gL$ and $\gL'$.

We present several applications of this result. Firstly, by combining
it with~\cite{BMR06}, we deduce an efficient algorithm that determines
whether two neighboring cluster-tilted algebras of any Dynkin type
(given by their quivers) are related by a BB-tilt. Secondly and more
generally, building on their Gorenstein property~\cite{KellerReiten07},
we give a numerical criterion for the derived equivalence
via BB-mutation of neighboring 2-CY-tilted algebras, stated only in
terms of their Cartan matrices, provided that they are invertible over
$\bQ$.

Another interesting class of algebras for which mutation can be related
to the quiver mutation of~\cite{FominZelevinsky02} consists of the
algebras of global dimension at most $2$.
Following~\cite{ABS08,Keller09}, we associate to such an algebra $A$ a
quiver $\wt{Q}_A$ called the extended quiver, and show that when it has
no loops and 2-cycles, then the extended quiver of a mutation
$\mu^-_k(A)$ or $\mu^+_k(A)$ that has global dimension at most $2$ is
obtained from $\wt{Q}_A$ by mutation at $k$.

We also provide an interpretation in terms of the generalized cluster
category $\cC_A$ of $A$ introduced in~\cite{Amiot08}. Namely, when $\cC_A$ is
$\Hom$-finite, then the image in $\cC_A$ of any of the complexes
$T^-_k$ and $T^+_k$  equals the mutation at $k$ of the canonical
cluster-tilting object in $\cC_A$, if the corresponding complex is
tilting with endomorphism algebra of global dimension at most $2$.

\section{Preliminaries}
\label{sec:pre}

Let $K$ be a field and let $A=KQ/I$ be a finite dimensional algebra
over $K$ given as a quotient of the path algebra of a finite quiver $Q$
by an admissible ideal $I$. For an arrow $\alpha$ in $Q$, let
$s(\alpha)$ and $t(\alpha)$ denote its start and end vertices. Our
convention is that two arrows $\alpha$, $\beta$ can be composed (as
$\alpha \beta$) if $t(\alpha) = s(\beta)$.

Let $\modf A$ be the category of finite dimensional right $A$-modules
and denote by $\dA$ its bounded derived category. Let $D = \Hom_k(-,k)$
be the duality on $k$-vector spaces. For a vertex $i$ of $Q$, denote by
$S_i$ the simple module corresponding to $i$ and by $P_i$ its
projective cover, which is spanned by all (non-zero) paths starting at
$i$. Thus an arrow $\alpha$ gives rise to a map $P_{t(\alpha)}
\xrightarrow{\alpha} P_{s(\alpha)}$.

Throughout this section, we fix a vertex $k$ of $Q$ which has no loops,
that is, $\Ext^1_A(S_k, S_k) = 0$.

\subsection{BB-tilting modules}
\label{ssec:BBtilt}

Let $\tau$ denote the Auslander-Reiten translation in $\modf A$.

\begin{defn}
\label{df:BBtilt}
We say that the BB-tilting module is \emph{defined} at the vertex $k$
if the $A$-module
\[
T^{\BB}_k = \tau^{-1} S_k \oplus \bigl( \bigoplus_{i \neq k} P_i \bigr)
\]
is a tilting module of projective dimension at most $1$. In this case,
$T^{\BB}_k$ is called the \emph{BB-tilting module associated with $k$}.
\end{defn}

It is shown in~\cite[Theorem~2.5]{TachikawaWakamatsu86} that
the condition in Definition~\ref{df:BBtilt} is
equivalent to the conditions that $\Hom_A(DA, S_k) = 0$ and
$\Ext^1_A(S_k, S_k) = 0$. In particular, $S_k$ is not injective,
meaning that $k$ is not a source. See also~\cite{BrennerButler80} for
the original construction and the survey article~\cite[\S2.8]{Assem90}.

\subsection{Perverse equivalence at a vertex}

We describe the construction of~\cite[\S2.6.2]{Rouquier06} (and
its dual) for the filtration $\phi \subset \{k\} \subset \{1, \dots,
n\}$ with perversities differing by $1$.

Let $k$ be a vertex of $Q$ without loops, and consider the maps
\begin{align*}
P_k \xrightarrow{f}
\bigoplus_{\alpha \,:\, t(\alpha) = k} P_{s(\alpha)}
&, &
\bigoplus_{\beta \,:\, s(\beta) = k} P_{t(\beta)}
\xrightarrow{g} P_k
\end{align*}
where $f=\bigoplus \alpha$ is induced by the arrows $\alpha$ ending at
$k$ and $g=\bigoplus \beta$ is induced by the arrows $\beta$ starting
at $k$. Note that $f$ is a minimal left $\add(A/P_k)$-approximation of
$P_k$ and $g$ is a minimal right $\add(A/P_k)$-approximation of $P_k$
(these notions are defined later in Section~\ref{sec:endo}).

Let $L_k$ be the cone of $f$ and $R_k$ be the cone of $g$ shifted one
place to the right. In a shortened notation, these are the complexes
\begin{align*}
L_k = P_k \xrightarrow{f} \bigoplus_{j \to k} P_j &, &
R_k = \bigoplus_{k \to j} P_j \xrightarrow{g} P_k
\end{align*}
where $P_k$ is in degree $-1$ in $L_k$ and in degree $1$ in $R_k$.

Define complexes $T^{-}_k$ and $T^{+}_k$ of projective $A$-modules by
\begin{align} \label{e:pTilt}
T^{-}_k = L_k \oplus \bigl( \bigoplus_{i \neq k} P_i \bigr) &, &
T^{+}_k = R_k \oplus \bigl( \bigoplus_{i \neq k} P_i \bigr) .
\end{align}
In other words, $T^{-}_k$ and $T^{+}_k$ are obtained by replacing the
summand $P_k$ in the (trivial) tilting complex $A$ by the complex $L_k$
or $R_k$, respectively. When we want to stress the dependency of these
complexes on the algebra $A$, we shall use the notation $T^{-}_k(A)$
and $T^{+}_k(A)$.

\section{BB-tilting vs.\ perverse equivalence at a vertex}
\label{sec:tilt}

We keep the assumptions of the previous section.
In particular, $k$ is a vertex without loops.

\subsection{Conditions for tilting}

In general, $T^{-}_k$ and $T^{+}_k$ need not be tilting complexes.
Therefore we start by giving the necessary and sufficient conditions on
$T^{-}_k$ and $T^{+}_k$ to be tilting, see also~\cite{HKM02},
\cite{Vitoria09} and~\cite[Theorem~4.1]{IyamaReiten08}.

Observe that the summands of each of the complexes $T^{-}_k$ and
$T^{+}_k$ always generate $\per A$, the triangulated subcategory of
$\dA$ consisting of the perfect complexes (that is, bounded complexes
of finitely generated projectives). This follows from the two short
exact sequences of complexes
\begin{align*}
0 \to \bigoplus_{j \to k} P_j \to L_k \to P_k[1] \to 0 &, & 0 \to
P_k[-1] \to R_k \to \bigoplus_{k \to j} P_j \to 0 ,
\end{align*}
recalling that none of the summands $P_j$ equals $P_k$.

\begin{prop} \label{p:pTiltCond}
Let $T^{-}_k$, $T^{+}_k$ be the complexes defined in~\eqref{e:pTilt}.
Then:
\begin{enumerate}
\renewcommand{\theenumi}{\alph{enumi}}
\item
\label{it:TLcpx}
$T^{-}_k$ is a tilting complex if and only if
\[
\Hom_{\dA}(P_i, L_k[-1]) = 0
\]
for all vertices $i \neq k$.

\item
\label{it:TLmod}
$T^{-}_k$ is isomorphic in $\dA$ to a tilting module if and only if
\[
\Hom_{\dA}(P_i, L_k[-1]) = 0
\]
for all vertices $i$.

\item
\label{it:TRcpx}
$T^{+}_k$ is a tilting complex if and only if
\[
\Hom_{\dA}(R_k[1], P_i) = 0
\]
for all vertices $i \neq k$.

\item
\label{it:TRmod}
$T^{+}_k$ is never isomorphic in $\dA$ to a tilting module.
\end{enumerate}
\end{prop}
\begin{proof}
The conditions in~(\ref{it:TLcpx}) and~(\ref{it:TRcpx}) are obviously
necessary. We show now claim~(\ref{it:TLcpx}). The proof
of~(\ref{it:TRcpx}) is similar.

Since $T^{-}_k$ is concentrated in only two consecutive degrees
it is enough to verify that
\[
\Hom_{\dA}(T_k^{-}, T_k^{-}[-1]) = 0 = \Hom_{\dA}(T_k^{-}, T_k^{-}[1]) .
\]
As we are dealing with complexes whose terms are projective, morphisms
in the derived category can be computed as morphisms in the homotopy
category of complexes. Obviously, the morphism spaces
\begin{align*}
\Hom(P_i, P_{i'}[1]) &,& \Hom(P_i, P_{i'}[-1]) &,& \Hom(P_i, L_k[1])
&,& \Hom(L_k, P_i[-1])
\end{align*}
vanish for all $i,i' \neq k$.

By assumption, $\Hom_{\dA}(P_i, L_k[-1]) = 0$ for every $i \neq k$. It
follows that $\Hom_{\dA}(L_k, L_k[-1]) = 0$ as well, by considering the
right square of a commutative diagram
\[
\xymatrix{
P_k \ar[r]^{f} \ar[d] & \bigoplus P_j \ar[r] \ar[d] & 0 \ar[d] \\
0 \ar[r] & P_k \ar[r]^{f} & \bigoplus P_j
}
\]
and using the assumption for each of the middle vertical maps $P_j \to P_k$,
recalling that none of the $P_j$ equals $P_k$.

In addition, $\Hom(L_k, P_i[1])=0$ for any $i \neq k$, as every path from
$i$ to $k$ (which corresponds to a map $P_k \to P_i$) factorizes
through one of the arrows ending at $k$, so that one can always define the
diagonal dotted homotopy in the diagram
\[
\xymatrix{P_k \ar[r]^{f} \ar[d] & \bigoplus P_j \ar[d] \ar@{.>}[dl] \\
P_i \ar[r] & 0
}
\]
This also shows that $\Hom(L_k, L_k[1]) = 0$, by considering the right
square in a commutative diagram
\[
\xymatrix{
0 \ar[r] \ar[d] & P_k \ar[r]^{f} \ar[d] & \bigoplus P_j \ar[d]
\ar@{.>}[dl] \\
P_k \ar[r]^{f} & \bigoplus P_j \ar[r] & 0
}
\]
and applying the above argument for each of the middle vertical maps
$P_k \to P_j$, recalling that none of the $P_j$ equals $P_k$. The proof
of~(\ref{it:TLcpx}) is thus complete.

For the proof of~(\ref{it:TLmod}), note that $\h^{-1}(T^{-}_k) \simeq
\Hom_{\dA}(A, L_k[-1])$.
Finally, for~(\ref{it:TRmod}) observe that the map $\bigoplus_{k \to j} P_j
\to P_k$ can never be surjective, hence $\h^1(T^{+}_k) \neq 0$.
\end{proof}

One can restate the conditions for $T^{-}_k$ to be a tilting complex in
terms of the kernel of the map $f$ whose cone is $L_k$.

\begin{cor} \label{c:TLkerf}
Let $T^{-}_k$ be the complex of~\eqref{e:pTilt}. Then:
\begin{enumerate}
\renewcommand{\theenumi}{\alph{enumi}}
\item \label{it:TLkerf}
$T^{-}_k$ is a tilting complex if and only if $\ker f \simeq S^m_k$ for
some $m \geq 0$.

\item \label{it:TLkerf0}
$T^{-}_k$ is isomorphic in $\dA$ to a tilting module if and only if $f$
is a monomorphism.
\end{enumerate}
\end{cor}
\begin{proof}
The claims follow from the isomorphisms
\[
\Hom_{\dA}(P_i, L_k[-1]) \simeq \Hom_A(P_i, \ker f)
\]
for any vertex $i$ and the corresponding claims in
Proposition~\ref{p:pTiltCond}.
Note that our assumption that are no loops at $k$ implies that any
module whose composition factors consist only of $S_k$ is of the
form $S_k^m$ for some $m \geq 0$.
\end{proof}

The conditions in Proposition~\ref{p:pTiltCond} can also be
conveniently rephrased in terms of (non-vanishing of) paths, as
follows. This is useful in practical calculations.

\begin{prop} \label{p:Tcomb}
Let $k$ be a vertex without loops and $T^{-}_k$, $T^{+}_k$ as
in~\eqref{e:pTilt}.
\begin{enumerate}
\renewcommand{\theenumi}{\alph{enumi}}
\item
\label{it:TLcomb}
$T^{-}_k$ is a tilting complex if and only if for any non-zero linear
combination $\sum a_r p_r$ of paths starting at $k$ and ending at some
vertex $i \neq k$, there exists at least one arrow $\alpha$ ending at
$k$ such that the composition $\sum a_r \alpha p_r$ is not zero.

\item
$T^{-}_k$ is isomorphic in $\dA$ to a tilting module if and only if the
condition of~(\ref{it:TLcomb}) holds for any vertex $i$ (including
$k$).

\item
$T^{+}_k$ is a tilting complex if and only if for any non-zero linear
combination of paths $\sum a_r p_r$ starting at some vertex $i \neq k$
and ending at $k$, there exists at least one arrow $\beta$ starting at
$k$ such that the composition $\sum a_r p_r \beta$ is not zero.
\end{enumerate}
\end{prop}
\begin{proof}
We show only~(\ref{it:TLcomb}), as the proof of the other assertions is
similar. By considering the diagram below,
\[
\xymatrix{
P_i \ar[r] \ar[d] & 0 \ar[d] \\
P_k \ar[r]^(0.4){f} & {\bigoplus_{j \to k} P_j}
}
\]
we see that the condition $\Hom_{\dA}(P_i, L_k[-1])=0$ is equivalent to
the condition that for any non-zero morphism $\rho : P_i \to P_k$,
there exists an arrow $\alpha : P_k \to P_j$ such that the composition
of $\rho$ with $\alpha$ gives a non-zero map $P_i \to P_j$.
\end{proof}

\begin{remark}
From the proposition we immediately see that if $k$ is a sink, then
$T^{-}_k$ is always a tilting complex while $T^{+}_k$ is never one,
whereas if $k$ is a source, then $T^{-}_k$ is never a tilting complex
while $T^{+}_k$ is always one.
\end{remark}

Let $Q^{op}$ be the opposite quiver of $Q$. Namely, it has the same set
of vertices as $Q$, with an (opposite) arrow $\alpha^{*} : j \to i$ for
any arrow $\alpha : i \to j$ of $Q$. If $A = KQ/I$, then the opposite
algebra $A^{op}$ can be written as $A^{op} = Q^{op}/I^{op}$ where
$I^{op}$ is generated by the paths opposite to those generating $I$.
The simple and indecomposable projective $A^{op}$-modules corresponding
to a vertex $i$ of $Q^{op}$ are then $D(S_i)$ and
$P^{*}_i = \Hom_A(P_i, A)$, respectively.

From the criteria in Proposition~\ref{p:Tcomb}, we immediately deduce
the following.

\begin{cor}
$T^{-}_k(A)$ is a tilting complex for $A$ if and only if
$T^{+}_k(A^{op})$ is a tilting complex for $A^{op}$.
\end{cor}

We conclude this section by giving another simple criterion for the
BB-tilting module to be defined at the vertex $k$. It will be used in
Section~\ref{ssec:numerical}.

\begin{lemma} \label{l:BBpres}
Let $k$ be a vertex without loops and let $\dots \to R^{-1} \to R^0 \to
0 \to \dots$ be a minimal projective resolution of $DA$. Then
$T^{\BB}_k$ is defined if and only if $P_k$ does not occur as summand
in $R^0$.
\end{lemma}
\begin{proof}
Recall that $T^{\BB}_k$ is defined if and only if $\Hom_A(DA, S_k)$
vanishes. The claim now follows from the fact that for any module $M$,
$\Hom_A(M, S_k) \simeq \Hom_A(P(M), S_k)$ where $P(M)$ is the
projective cover of $M$.
\end{proof}

\subsection{BB-tilt as perverse equivalence at a vertex and vice versa}

\begin{lemma}
The sequence
\begin{equation} \label{e:tSk}
P_k \xrightarrow{f} \bigoplus_{j \to k} P_j \to \tau^{-1} S_k \to 0
\end{equation}
is a projective presentation of $\tau^{-1} S_k$.
\end{lemma}
\begin{proof}
Since $\tau^{-1} = \Tr D$, we start by writing a projective
presentation of the simple $A^{op}$-module $D(S_k)$ as
\[
\bigoplus_{k \to j} P^{*}_j \xrightarrow{f^*} P^{*}_k \to D(S_k) \to 0
\]
where the sum goes over the arrows $\alpha^{*} : k \to j$ in $Q^{op}$,
and apply $\Tr$ to get~\eqref{e:tSk}.
\end{proof}

The relations between BB-tilting and perverse equivalences at a vertex
are summarized in the following proposition.

\begin{prop} \label{p:BBperv}
Let $k$ be a vertex (without loops) of $Q$.
\begin{enumerate}
\renewcommand{\theenumi}{\alph{enumi}}
\item \label{it:TBBp}
If the $BB$-tilting module is defined at $k$, then $T^{-}_k$ is a
tilting complex and $T^{\BB}_k \simeq T^{-}_k$ in $\dA$.

\item \label{it:TpBB}
Conversely, if $T^{-}_k$ is isomorphic in $\dA$ to a tilting module,
then the $BB$-tilting module is defined at $k$ and $T^{-}_k \simeq
T^{\BB}_k$ in $\dA$.
\end{enumerate}
\end{prop}
\begin{proof}
\begin{enumerate}
\renewcommand{\theenumi}{\alph{enumi}}
\item
The assumption implies that $\pd_A \tau^{-1} S_k \leq 1$. It follows
that the projective presentation of~\eqref{e:tSk} is actually a
resolution. In other words, $\tau^{-1} S_k \simeq L_k$ in $\dA$. Hence
$T^{\BB}_k \simeq T^{-}_k$.

\item
By Corollary~\ref{c:TLkerf}, $\ker f = 0$, hence~\eqref{e:tSk} is a
projective resolution of $\tau^{-1} S_k$. Therefore
$\pd_A \tau^{-1}S_k \leq 1$ and $\tau^{-1} S_k \simeq L_k$ in $\dA$.
Hence $T^{\BB}_k \simeq T^{-}_k$ is a tilting module of projective
dimension at most $1$.
\end{enumerate}
\end{proof}

Under certain additional assumptions, any perverse equivalence given by
a tilting complex $T^{-}_k$ can be regarded as a BB-tilting module.

\begin{lemma} \label{l:SradP}
Assume that $S_k$ is not a submodule of the radical of $P_k$. If
$T^{-}_k$ is a tilting complex, then the BB-tilting module is defined
at $k$ and $T^{-}_k \simeq T^{\BB}_k$ in $\dA$.
\end{lemma}
\begin{proof}
By Corollary~\ref{c:TLkerf}(\ref{it:TLkerf}), $\ker f \simeq S_k^m$ for
some $m \geq 0$. On the other hand, $\ker f \subseteq \rad P_k$.
Therefore $\ker f = 0$ and the result follows from
Corollary~\ref{c:TLkerf}(\ref{it:TLkerf0}) and
Proposition~\ref{p:BBperv}(\ref{it:TpBB}).
\end{proof}

\begin{remark}
When $\End_A(P_k) \simeq K$, the assumption of the lemma holds. This
happens, in particular, when the quiver $Q$ is acyclic (i.e.\ $A$ is
\emph{triangular}), or when the algebra $A$ is \emph{schurian}, that
is, $\dim_K \Hom_A(P_i, P_{i'}) \leq 1$ for any two vertices $i, i'$ of
$Q$.
\end{remark}

\begin{remark}
When $A$ is a symmetric algebra which is not semi-simple, $T^{-}_k$ is
always a tilting complex~\cite[\S2.6]{Rouquier06}, but the assumption
of Lemma~\ref{l:SradP} does not hold. Observe that $T^{\BB}_k$ is never
defined since $A$ has no non-trivial tilting modules.
\end{remark}

\begin{lemma}
Assume that $\pd_A \tau^{-1} S_k \leq 2$. If $T^{-}_k$ is a tilting
complex, then the BB-tilting module is defined at $k$ and $T^{-}_k
\simeq T^{\BB}_k$ in $\dA$.
\end{lemma}
\begin{proof}
By Corollary~\ref{c:TLkerf}(\ref{it:TLkerf}), the
sequence~\eqref{e:tSk} yields an exact sequence
\[
0 \to S_k^m \to P_k \xrightarrow{f} \bigoplus_{j \to k} P_j \to
\tau^{-1} S_k \to 0
\]
for some $m \geq 0$. If $m>0$, then by our assumption, $S_k$ must be
projective, hence $P_k = S_k$, so that $f$ is a monomorphism, a
contradiction. Thus $m=0$ and the result follows.
\end{proof}

\section{Mutations of algebras}
\label{sec:algmut}

\subsection{Operations}
We assume now that the field $K$ is algebraically closed. Let $A =
KQ/I$ be a finite dimensional $K$-algebra given as a quiver $Q$ with
relations.

\begin{defn}
Let $k$ be a vertex of $Q$ without loops.
\begin{enumerate}
\renewcommand{\theenumi}{\alph{enumi}}
\item
We say that the negative mutation is \emph{defined} at the vertex $k$
if $T^{-}_k(A)$ is a tilting complex. In this case, we call
\[
\mu^{-}_k(A) = \End_{\dA} T^{-}_k(A)
\]
the \emph{negative mutation of $A$ at the vertex $k$}.

\item
We say that the positive mutation is \emph{defined} at the vertex $k$
if $T^{+}_k(A)$ is a tilting complex. In this case,
\[
\mu^{+}_k(A) = \End_{\dA} T^{+}_k(A)
\]
is called the \emph{positive mutation of $A$ at the vertex $k$}.

\item
We say that the BB-mutation is \emph{defined} at the vertex $k$ if
$T^{\BB}_k(A)$ is defined. In this case, we call
\[
\mu^{\BB}_k(A) = \End_{A} T^{\BB}_k(A)
\]
the \emph{BB-mutation of $A$ at the vertex $k$}.
\end{enumerate}
\end{defn}

Obviously, when $\mu^{\BB}_k(A)$ is defined, so is $\mu^-_k(A)$, and
moreover $\mu^{\BB}_k(A) \simeq \mu^-_k(A)$. The following proposition
justifies the name ``mutations'' for these operations.

\begin{prop}
Let $k$ be a vertex of $Q$ without loops.
\begin{enumerate}
\renewcommand{\theenumi}{\alph{enumi}}
\item
If $\mu^{-}_k(A)$ is defined, then $\mu^{+}_k(\mu^{-}_k(A))$ is defined
and isomorphic to $A$.

\item
If $\mu^{+}_k(A)$ is defined, then $\mu^{-}_k(\mu^{+}_k(A))$ is defined
and isomorphic to $A$.

\item
If $\mu^{\BB}_k(A)$ is defined, then
$\mu^{\BB}_k\bigl(\left(\mu^{\BB}_k(A)\right)^{op}\bigr)$ is defined
and isomorphic to $A^{op}$.
\end{enumerate}
\end{prop}

Given a vertex $k$, precisely one of the following situations can
occur:
\begin{itemize}
\item
Both $\mu^{-}_k(A)$ and $\mu^{+}_k(A)$ are defined; or
\item
One of them is defined; or
\item
None of them is defined.
\end{itemize}
When $k$ is a sink or a source in $Q$, then exactly one of the
mutations is defined, namely, $\mu^-_k(A)$ for a sink $k$ and
$\mu^+_k(A)$ for a source $k$.

\begin{example} \label{ex:algmut}
Let $A$ be the path algebra of the quiver
\[
\xymatrix{{\bullet_1} \ar[r] & {\bullet_2} \ar[r] & {\bullet_3}} .
\]

Using Proposition~\ref{p:Tcomb}, we see that at the vertex $1$, only
$\mu^+_1(A)$ is defined, whereas at the vertex $3$, only $\mu^-_3(A)$
is defined. At the vertex $2$, both $\mu^-_2(A)$ and $\mu^+_2(A)$ are
defined, and are given by the following quivers with zero relations:
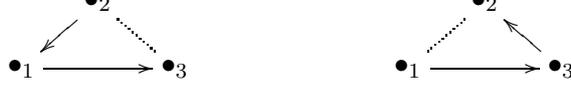
\begin{align*}
\xymatrix@=1pc{
& {\bullet_2} \ar[ld] \\
{\bullet_1} \ar[rr] & & {\bullet_3} \ar@{.}[lu]
}
& &
\xymatrix@=1pc{
& {\bullet_2} \ar@{.}[ld] \\
{\bullet_1} \ar[rr] & & {\bullet_3} \ar[lu]
}
\end{align*}

Let $A'=\mu^-_2(A)$ be as in the left picture. Then at the vertex $1$,
none of $\mu^-_1(A')$ and $\mu^+_1(A')$ is defined.
\end{example}

\subsection{K-theoretical interpretation}

Let $Q$ be a quiver with $n$ vertices. For a vertex $1 \leq k \leq n$,
define the $n \times n$ matrices $r^{-}_k = r^{-}_k(Q)$ and $r^{+}_k =
r^{+}_k(Q)$ by
\begin{align*}
\left(r^{-}_k\right)_{ij} &=
\begin{cases}
-\delta_{ij} + \left| \{\text{arrows $j \to k$ in $Q$}\} \right| &
\text{if $i=k$,} \\
\delta_{ij} & \text{otherwise,}
\end{cases}
\intertext{and}
\left(r^{+}_k\right)_{ij} &=
\begin{cases}
-\delta_{ij} + \left| \{\text{arrows $k \to j$ in $Q$}\} \right| &
\text{if $i=k$,} \\
\delta_{ij} & \text{otherwise.}
\end{cases}
\end{align*}

As already observed in~\cite[Lemma~7.1]{GLS06}, these matrices are
closely related to the Fomin-Zelevinsky matrix
mutation~\cite{FominZelevinsky02} of the skew-symmetric matrix
corresponding to $Q$. Namely, recall that for a quiver $Q$ there is an
associated skew-symmetric matrix $b_Q$ defined by
\begin{equation} \label{e:BfromQ}
(b_Q)_{ij} = \bigl| \{\text{arrows $j \to i$} \} \bigr|
           - \bigl| \{\text{arrows $i \to j$} \} \bigr|,
\end{equation}
and one can recover $Q$ from $b_Q$ as long as it has no loops or
$2$-cycles.

Fomin and Zelevinsky have defined the mutation $\mu_k(b_Q)$, which is a
again a skew-symmetric matrix, for any vertex $k$. Then $\mu_k(b_Q)$ is
obtained from $b_Q$ by viewing the latter as a matrix of a bilinear
form and applying a change-of-basis transformation given by the matrix
$r^-_k$ or $r^+_k$. More precisely, the following lemma is a
reformulation of~\cite[Lemma~7.1]{GLS06}.

\begin{lemma} \label{l:rmut}
Assume that $Q$ has no loops and no $2$-cycles. Then for any vertex
$k$, we have
\[
\mu_k(b_Q) = \left(r^-_k\right)^T b_Q r^-_k = \left(r^+_k \right)^T b_Q
r^+_k .
\]
\end{lemma}

Now we relate the matrices $r^-_k$ and $r^+_k$ to mutations of algebras
via the notion of the Euler form. Let $A=KQ/I$ be a finite dimensional
$K$-algebra given as a quiver with relations and let $n$ be the number
of vertices of $Q$. The \emph{Cartan matrix} of $A$ is the $n \times n$
integral matrix $C_A$ whose entries are
\[
(C_A)_{ij} = \dim_K \Hom_A(P_i, P_j)
\]
for $1 \leq i, j \leq n$.

Recall that the Grothendieck group $K_0(\per A)$ is free abelian on the
generators $[P_1], \dots, [P_n]$, and the expression
\begin{equation} \label{e:EulerForm}
\langle X, Y \rangle_A = \sum_{r \in \bZ} (-1)^r \dim_K \Hom_{\dA}(X,
Y[r])
\end{equation}
is well defined for any $X, Y \in \per A$ and induces a bilinear form
on $K_0(\per A)$, known as the \emph{Euler form}, whose matrix with
respect to the basis of projectives is $C_A$.

The following lemma is briefly mentioned in~\cite{BHL09}.

\begin{lemma} \label{l:TCartan}
Let $T$ be a (basic) tilting complex in $\dA$ with endomorphism algebra
$A' = \End_{\dA}(T)$ and let $T_1, \dots, T_n$ be the indecomposable
summands of $T$. Then the Cartan matrix of $A'$ is given by $C_{A'} = r
C_A r^T$, where $r = (r_{ij})_{i,j=1}^n$ is the matrix defined by
\[
[T_i] = \bigoplus_{j=1}^n r_{ij} [P_j]
\]
(that is, its $i$-th row is the class of the summand $T_i$ in $K_0(\per
A)$ written in the basis $[P_1], \dots, [P_n]$).
\end{lemma}
\begin{proof}
As the indecomposable projectives of $A'$ correspond to the
indecomposable summands of $T$, we have
\[
\left( C_{A'} \right)_{ij} = \dim_K \Hom_{\dA}(T_i, T_j) = \langle
[T_i], [T_j] \rangle_A
\]
since $T$ is a tilting complex. The result now follows.
\end{proof}

\begin{prop} \label{p:rCartan}
Let $A=KQ/I$ be a finite dimensional $K$-algebra.
\begin{enumerate}
\renewcommand{\theenumi}{\alph{enumi}}
\item
If $\mu^{-}_k(A)$ is defined, then $C_{\mu^{-}_k(A)} = r^{-}_k \cdot
C_A \cdot \left( r^{-}_k \right)^T$.

\item
If $\mu^{+}_k(A)$ is defined, then $C_{\mu^{+}_k(A)} = r^{+}_k \cdot
C_A \cdot \left( r^{+}_k \right)^T$.
\end{enumerate}
\end{prop}
\begin{proof}
Use Lemma~\ref{l:TCartan} and the definition of $T^{-}_k(A)$ and
$T^{+}_k(A)$ in~\eqref{e:pTilt}.
\end{proof}

When the algebra $A$ has finite global dimension, $\per A = \dA$ and
the classes $[S_1], \dots, [S_n]$ form a basis of $K_0(\per A)$. The
matrix of the Euler form~\eqref{e:EulerForm} with respect to that basis
is then given by $c_A = C_A^{-T}$, that is, the inverse of the
transpose of $C_A$.

\begin{cor} \label{c:rEuler}
Assume that $A$ has finite global dimension.
\begin{enumerate}
\renewcommand{\theenumi}{\alph{enumi}}
\item
If $\mu^{-}_k(A)$ is defined, then
$c_{\mu^{-}_k(A)} = \left(r^{-}_k \right)^T \cdot c_A \cdot r^{-}_k$.

\item
If $\mu^{+}_k(A)$ is defined, then
$c_{\mu^{+}_k(A)} = \left(r^{+}_k \right)^T \cdot c_A \cdot r^{+}_k$.
\end{enumerate}
\end{cor}
\begin{proof}
By Proposition~\ref{p:rCartan} we have
\begin{align*}
c_{\mu^-_k(A)} &= \left( C_{\mu^-_k(A)} \right)^{-T} = \left( r^-_k C_A
\left( r^-_k \right)^T \right)^{-T} = \left( r^-_k \right)^{-T} c_A
\left( r^-_k \right)^{-1} \\
&= \left( r^-_k \right)^T c_A r^-_k
\end{align*}
where for the last equality we used the fact that $\left( r^-_k
\right)^2$ is the identity matrix.
\end{proof}

\section{Mutations of endomorphism algebras}
\label{sec:endo}

Mutations of algebras arise naturally when considering endomorphism
algebras of objects related by approximation sequences in additive
categories. This has been essentially observed by Hu and Xi
in~\cite{HuXi08}, where they introduced the notion of almost
$\cD$-split sequences. We recall their result, formulating it in a form
which will be convenient for our applications.

Let $\cC$ be a category, $\cD$ a full subcategory and $X$ an object of
$\cC$. A morphism $f : X \to D$ is called a \emph{left $\cD$-approximation}
if $D \in \cD$ and any morphism $f' : X \to D'$ with $D' \in \cD$ can
be completed to a commutative diagram as in the left picture.
\begin{align*}
\xymatrix{ X \ar[rd]_{f'} \ar[r]^{f} & D \ar@{.>}[d] \\ & D' } &&
\xymatrix{ D \ar[r]^{g} & X  \\ D' \ar[ru]_{g'} \ar@{.>}[u] }
\end{align*}
$f$ is a \emph{minimal left $\cD$-approximation} if furthermore, when
considering the left diagram with $f'=f$, any dotted arrow making it
commutative is an automorphism of $D$. The notions of a \emph{right
$\cD$-approximation} and a \emph{minimal right $\cD$-approximation} are
defined similarly using the right diagram, see~\cite{AuslanderSmalo80}.

Let $\cC$ be an additive category with split idempotents. For an object
$M \in \cC$, denote by $\add M$ the full subcategory consisting of
finite direct sums of direct summands of $M$. It is equivalent to the
additive category of finitely generated projective
$\End_{\cC}(M)$-modules, via the functor $\Hom_{\cC}(M,-)$.

\begin{prop}[Lemma~3.4 of~\protect{\cite{HuXi08}}] \label{p:EndoBB}
Let $M \in \cC$ and let
\[
X \xrightarrow{f} B \xrightarrow{g} X'
\]
be a sequence of morphisms satisfying the following conditions:
\begin{enumerate}
\renewcommand{\theenumi}{\roman{enumi}}
\item
$f : X \to B$ is a left $\add M$-approximation of $X$
and $g : B \to X'$ is a right $\add M$-approximation of $X'$,

\item
The induced sequences
\begin{align*}
& 0 \to \Hom_{\cC}(U, X) \xrightarrow{f_*} \Hom_{\cC}(U, B)
\xrightarrow{g_*} \Hom_{\cC}(U, X') \\
& 0 \to \Hom_{\cC}(X', U') \xrightarrow{g^*} \Hom_{\cC}(B, U')
\xrightarrow{f^*} \Hom_{\cC}(X, U')
\end{align*}
are exact, where $U = X \oplus M$ and $U' = X' \oplus M$.
\end{enumerate}
Then the rings $\gL = \End_{\cC}(U)$ and $\gL' = \End_{\cC}(U')$ are
derived equivalent.
\end{prop}

Moreover, examining the proof in~\cite{HuXi08}, we see that the
following complex of projective $\gL$-modules
\begin{equation} \label{e:Endcpx}
\left(\Hom_{\cC}(U,X) \xrightarrow{f_*} \Hom_{\cC}(U,B) \right) \oplus
\Hom_{\cC}(U, M) ,
\end{equation}
where $\Hom_{\cC}(U,X)$ is in degree $-1$ and the other terms are in
degree $0$, is a tilting complex whose endomorphism ring is isomorphic
to $\gL'$. The exactness of the first sequence implies that this
complex is quasi-isomorphic to a tilting module.

The similarity of~\eqref{e:Endcpx} and~\eqref{e:pTilt} suggests that
one can replace some of the conditions on $f_*$ and $g^*$ by intrinsic
conditions expressed only in terms of the algebras $\gL$ and $\gL'$.

\begin{thm} \label{t:Endoder}
Let $K$ be an algebraically closed field and $\cC$ be a $K$-linear,
$\Hom$-finite, additive category with split idempotents. Let $U_1, U_2,
\dots, U_{n-1}$ be non-isomorphic indecomposable objects of $\cC$ and
set $M = U_1 \oplus \dots \oplus U_{n-1}$.

Let $U_n$ and $U'_n$ be indecomposable objects of $\cC$ not in $\add
M$. Set $U = M \oplus U_n$ and $U' = M \oplus U'_n$ and consider the
endomorphism algebras
\[
\gL = \End_{\cC}(U) \quad \text{ and } \quad \gL' = \End_{\cC}(U') .
\]
Number the vertices of their quivers in accordance with the numbering
of the summands $U_i$. Assume that:
\begin{enumerate}
\renewcommand{\theenumi}{\roman{enumi}}
\item \label{it:approx}
There exist morphisms
\[
f : U_n \to B \quad \text{ and } \quad g : B \to U'_n
\]
such that $f$ is a minimal left $\add M$-approximation and $g$ is a
minimal right $\add M$-approximation;

\item \label{it:Homseqex}
The induced sequences
\begin{align*}
\Hom_{\cC}(U, U_n) \xrightarrow{f_*} &\Hom_{\cC}(U, B)
\xrightarrow{g_*} \Hom_{\cC}(U, U'_n) \\
\Hom_{\cC}(U'_n, U') \xrightarrow{g^*} &\Hom_{\cC}(B, U')
\xrightarrow{f^*} \Hom_{\cC}(U_n, U')
\end{align*}
are exact.
\end{enumerate}
Then the following are equivalent:
\begin{enumerate}
\renewcommand{\theenumi}{\alph{enumi}}
\item \label{it:Endoder:BB}
The $BB$-tilting modules $T^{\BB}_n(\gL)$ and $T^{\BB}_n(\gL'^{op})$
are defined (in particular, the quivers of $\gL$ and $\gL'$ have no
loops at the vertex $n$).

\item \label{it:Endoder:fg}
The induced maps $f_*$ and $g^*$ are monomorphisms.

\item \label{it:Endoder:der}
$\gL' \simeq \mu^{\BB}_n(\gL)$ (in particular, $\gL$ and $\gL'$ are
derived equivalent).
\end{enumerate}
\end{thm}
\begin{proof}
The indecomposable projectives in $\modf \gL$ are precisely the modules
$\Hom_{\cC}(U, U_i)$. The assumption that $f: U_n \to B$ is a minimal
left $\add M$-approximation implies that $f_* : \Hom_{\cC}(U, U_n) \to
\Hom_{\cC}(U, B)$ is a minimal left $\add \Hom_{\cC}(U,
M)$-approximation in $\modf \gL$, and therefore the complex
$T^-_n(\gL)$ can be written as
\[
T^-_n(\gL) \simeq \left(\Hom_{\cC}(U, U_n) \xrightarrow{f_*}
\Hom_{\cC}(U, B) \right) \oplus \Hom_{\cC}(U, M) .
\]
From Proposition~\ref{p:BBperv} and
Corollary~\ref{c:TLkerf}(\ref{it:TLkerf0}) we deduce that $f_*$ is a
monomorphism if and only if $T^{\BB}_n(\gL)$ is defined.

Similarly, the indecomposable projectives in $\modf \gL'^{op}$ are the
modules $\Hom_{\cC}(U_i, U')$ for $i<n$ together with $\Hom_{\cC}(U'_n,
U')$. The assumption that $g: B \to U'_n$ is a minimal right $\add
M$-approximation implies that $g^* : \Hom_{\cC}(U'_n, U') \to
\Hom_{\cC}(B, U')$ is a minimal left $\add \Hom_{\cC}(M,
U')$-approximation in $\modf \gL'^{op}$, and therefore the complex
$T^-_n(\gL'^{op})$ can be written as
\[
T^-_n(\gL'^{op}) \simeq \left(\Hom_{\cC}(U'_n, U') \xrightarrow{g^*}
\Hom_{\cC}(B, U') \right) \oplus \Hom_{\cC}(M, U') .
\]
A similar reasoning to the above shows that $g^*$ is a monomorphism if
and only if $T^{\BB}_n(\gL'^{op})$ is defined. This shows the
equivalence of (\ref{it:Endoder:BB}) and (\ref{it:Endoder:fg}).

Now (\ref{it:Endoder:der}) follows from (\ref{it:Endoder:fg}) by
Proposition~\ref{p:EndoBB} and the discussion afterwards. Finally, if
$\gL' \simeq \mu^{\BB}_n(\gL)$, then $\gL^{op} \simeq
\mu^{\BB}_n(\gL'^{op})$ and hence (\ref{it:Endoder:der}) implies
(\ref{it:Endoder:BB}).
\end{proof}

We see that if $\gL$ and $\gL'$ are the endomorphism algebras of two
objects related by replacing an indecomposable summand by another
through an approximation sequence, then the existence of the two
BB-mutations $\mu^{\BB}_n(\gL)$ and $\mu^{\BB}_n(\gL'^{op})$
automatically implies that they take the ``correct'' values, namely
$\gL'$ and $\gL^{op}$, respectively.

In the applications, the condition~(\ref{it:Homseqex}) in the theorem
is usually automatically satisfied. Indeed, consider an approximation
sequence
\begin{equation} \label{e:approx}
U_n \xrightarrow{f} B \xrightarrow{g} U'_n
\end{equation}
as in~(\ref{it:approx}) of the theorem.

\begin{cor} \label{c:Endoder}
If $\cC$ is exact and~\eqref{e:approx} is an exact sequence,
then~(\ref{it:Homseqex}) is satisfied. Moreover, in this case $f_*$ and
$g^*$ are monomorphisms so that always $\gL' \simeq \mu^{\BB}_n(\gL)$.
\end{cor}

\begin{remark}
If $\cC$ is triangulated and~\eqref{e:approx} is a triangle,
then~(\ref{it:Homseqex}) is satisfied. In this case, the BB-tilting
modules $T^{\BB}_n(\gL)$, $T^{\BB}_n(\gL'^{op})$ are not always defined
and the algebras $\gL$, $\gL'$ are not necessarily derived equivalent.
\end{remark}

\section{Mutations of endomorphism algebras in 2-CY categories}
\label{sec:2CY}

Cluster categories were introduced in~\cite{BMRRT06} as a categorical
framework for the cluster algebras of Fomin and
Zelevinsky~\cite{FominZelevinsky02}. Particular role is played by the
so-called cluster-tilting objects and their endomorphism algebras,
known as cluster-tilted algebras~\cite{BMR07}. Cluster categories are
2-Calabi-Yau triangulated categories, and the theory has been extended
to such
categories~\cite{Amiot08,BIRSc09,BIRSm08,IyamaYoshino08,KellerReiten07}
as well as to Frobenius categories which are stably
2-Calabi-Yau~\cite{BIRSc09,GLS06}.

The results of the preceding section can be applied in the study of
endomorphism algebras of cluster-tilting objects in these categories.
We begin by recalling the relevant notions.

Let $K$ be an algebraically closed field and $\cC$ a $K$-linear
$\Hom$-finite additive category with split idempotents which is either:
\begin{itemize}
\item
triangulated 2-Calabi-Yau (2-CY) category, that is, there exist functorial
isomorphisms
$\Hom_{\cC}(X, Y[2]) \simeq D\Hom_{\cC}(Y, X)$ for $X, Y \in \cC$; or

\item
a Frobenius category whose stable category, which is triangulated
by~\cite{Happel88}, is 2-CY.
\end{itemize}
When $\cC$ is triangulated, we denote $\Hom_{\cC}(X,Y[1])$ by
$\Ext^1_{\cC}(X,Y)$.

An object $U \in \cC$ is \emph{rigid} if $\Ext^1_{\cC}(U, U) = 0$. It
is a \emph{cluster-tilting object} if it is rigid and for any $X \in
\cC$, $\Ext^1_{\cC}(U, X) = 0$ implies that $X \in \add U$.
When $\cC$ is Frobenius stably 2-CY, any indecomposable
projective-injective is a summand of any cluster-tilting object.

The operation of \emph{mutation} of~\cite{FominZelevinsky02} is
categorified by mutation of cluster-tilting objects, see~
\cite{BMRRT06,IyamaYoshino08}. Let $U = U_1 \oplus U_2 \oplus \dots
\oplus U_n$ be a basic cluster-tilting object written as a sum of $n$
non-isomorphic indecomposables.
When $\cC$ is triangulated, set $m=n$. When $\cC$ is Frobenius,
assume that $U_{m+1}, \dots, U_n$ are all the projective-injective
summands of $U$. Then for any $1 \leq k \leq m$ there
exists a unique $U'_k$ non-isomorphic to $U_k$ such that $\mu_k(U) =
(U/U_k) \oplus U'_k$ is a cluster-tilting object (here, $U/U_k$ denotes
the sum of all $U_i$ for $i \neq k$). Moreover, there are so-called
\emph{exchange triangles}
\begin{align} \label{e:exchtri}
U_k \xrightarrow{f} B \xrightarrow{g} U'_k && U'_k \xrightarrow{f'} B'
\xrightarrow{g'} U_k
\end{align}
when $\cC$ is triangulated, and (exact) \emph{exchange sequences}
\begin{align} \label{e:exchseq}
0 \to U_k \xrightarrow{f} B \xrightarrow{g} U'_k \to 0 &&
0 \to U'_k \xrightarrow{f'} B' \xrightarrow{g'} U_k \to 0
\end{align}
when $\cC$ is Frobenius,
such that the maps $f, f'$ are minimal left $\add (U/U_k)$-approximations
and $g, g'$ are minimal right $\add (U/U_k)$-approximations.

\subsection{Endomorphism algebras in 2-CY Frobenius categories}

Let $K$ be an algebraically closed field and
$\cE$ a $K$-linear, $\Hom$-finite, Frobenius category which is stably
2-CY.

It is known by~\cite{Palu09} that the derived equivalence class of the
endomorphism algebra of a cluster-tilting object $U$ in $\cE$ does not
depend on $U$, see also~\cite[\S5]{Iyama07}. It is still of interest to
relate the endomorphism algebras of two cluster-tilting objects related
by mutation. This has been done in~\cite{GLS06} for Frobenius
categories arising from preprojective algebras. The following is a
generalization, essentially saying that the mutation operations of
algebras ``commute'' with mutation of cluster-tilting objects.

\begin{thm}
Let $U$ be a cluster-tilting object in $\cE$.
Then for any non projective-injective summand $U_k$ of $U$,
all the corresponding BB-, negative and positive mutations of
$\End_{\cC}(U)$ are defined and moreover,
\[
\End_{\cE}\bigl(\mu_k(U)\bigr) \simeq
\mu^{\BB}_k\bigl(\End_{\cE}(U)\bigr) =
\mu^{-}_k\bigl(\End_{\cE}(U)\bigr) \simeq
\mu^{+}_k\bigl(\End_{\cE}(U)\bigr) .
\]
\end{thm}
\begin{proof}
Using Corollary~\ref{c:Endoder} for the left exchange sequence
in~\eqref{e:exchseq}, we get that $\mu^{\BB}_k(\End_{\cC}(U))$
is defined and
\[
\End_{\cE}\bigl(\mu_k(U)\bigr) \simeq
\mu^{\BB}_k\bigl(\End_{\cE}(U)\bigr) =
\mu^{-}_k\bigl(\End_{\cE}(U)\bigr) .
\]
Using this for $\mu_k(U)$ instead of $U$, we get that
\[
\End_{\cE}(U) \simeq
\End_{\cE}\bigl(\mu_k(\mu_k(U))\bigr) \simeq
\mu^{\BB}_k\bigl(\End_{\cE}(\mu_k(U))\bigr) =
\mu^{-}_k\bigl(\End_{\cE}(\mu_k(U))\bigr) ,
\]
hence $\mu^+_k(\End_{\cE}(U))$ is also defined and isomorphic to
$\End_{\cC}(\mu_k(U))$.
\end{proof}

\subsection{Endomorphism algebras in 2-CY triangulated categories}

Let $K$ be an algebraically closed field and
$\cC$ a $K$-linear, $\Hom$-finite, triangulated 2-CY category with
split idempotents.

Let $U = U_1 \oplus U_2 \oplus \dots \oplus U_n$ a cluster-tilting
object.
According to~\cite{BMR07,KellerReiten07}, the neighboring 2-CY-tilted
algebras $\gL = \End_{\cC}(U)$ and $\gL' = \End_{\cC}(\mu_k(U))$ are
related by what is called in~\cite{Ringel07} \emph{near-Morita
equivalence}, namely if $S_k$ and $S'_k$ are the simple modules in
$\gL$ and $\gL'$ corresponding to the indecomposables $U_k$ and $U'_k$,
then
\[
\modf \gL / \langle \add S_k \rangle \simeq \modf \gL' / \langle \add
S_k' \rangle .
\]
This generalizes APR-tilting~\cite{APR79}, which corresponds to the
case where the vertex $k$ is a sink.

Another feature of APR-tilting is that the algebras related by an
APR-tilt are derived equivalent. However, it is easily seen (and well
known) that two nearly-Morita equivalent 2-CY-tilted algebras are in
general not derived equivalent.

\begin{example} \label{ex:A3}
The quivers of two cluster-tilted algebras of type $A_3$ are shown
below. One is obtained from the other by mutation at the vertex $2$.
\begin{align} \label{e:ctaA3}
\xymatrix@=1pc{
& {\bullet_2} \ar[rd] \\
{\bullet_1} \ar[ru] & & {\bullet_3} } & & \xymatrix@=1pc{
& {\bullet_2} \ar[ld] \\
{\bullet_1} \ar[rr] & & {\bullet_3} \ar[lu] }
\end{align}
The corresponding algebras, where in the right one the composition of
any pair of consecutive arrows is zero, are nearly Morita equivalent
but not derived equivalent.
\end{example}

Nevertheless, several derived equivalence classifications of
cluster-tilted algebras~\cite{Bastian09,BHL09,BuanVatne08} have
revealed that, at least in finite mutation type, the number of derived
equivalence classes is much smaller than the number of isomorphism
classes of algebras (which equals the number of quivers in the mutation
class). For example, in type $E_8$ there are $1574$ algebras but only
$15$ derived equivalence classes~\cite{BHL09}. Moreover, many of the
classifications rely on showing that sufficiently many pairs of near
Morita equivalent cluster-tilted algebras are also derived equivalent.

In an attempt to provide a conceptual explanation of these phenomena,
it is therefore interesting to formulate conditions that will guarantee
the derived equivalence of neighboring 2-CY-tilted algebras. A similar
problem for Jacobian algebras was studied in~\cite{Vitoria09}.

\begin{thm} \label{t:2CYtri}
Let $U$ be a cluster-tilting object in $\cC$ and let $\gL =
\End_{\cC}(U)$ and $\gL' = \End_{\cC}(\mu_k(U))$ be two neighboring
2-CY-tilted algebras.

Then $\gL' \simeq \mu^{\BB}_k(\gL)$ if and only if the BB-tilting
modules $T^{\BB}_k(\gL)$ and $T^{\BB}_k(\gL'^{op})$ are defined. In
particular, in this case $\gL$ and $\gL'$ are derived equivalent.
\end{thm}
\begin{proof}
Use the left exchange triangle in~\eqref{e:exchtri} and
Theorem~\ref{t:Endoder}.
\end{proof}

\begin{example}
Looking again at Example~\ref{ex:A3}, denote by $\gL$ and $\gL'$ the
left and right cluster-tilted algebras in~\eqref{e:ctaA3}. Then
$T^{\BB}_2(\gL)$ is defined and $\mu^{\BB}_2(\gL)$ is the algebra given
by the left quiver of Example~\ref{ex:algmut} (with zero relation),
hence $\gL' \not \simeq \mu^{\BB}_2(\gL)$. Observe that for any $1 \leq
i \leq 3$, none of the mutations $\mu^+_i(\gL')$, $\mu^-_i(\gL')$ is
defined.
\end{example}

Theorem~\ref{t:2CYtri} motivates the following definition, see
also~\cite{BHL09}.
\begin{defn}
Let $U$ be a cluster-tilting object of $\cC$. The mutation $\mu_k(U)$
is \emph{good} if $\End_{\cC}(\mu_k(U)) \simeq
\mu^{\BB}_k(\End_{\cC}(U))$. In this case, we say that the algebra
$\End_{\cC}(\mu_k(U))$ is obtained from $\End_{\cC}(U)$ by a \emph{good
mutation}.
\end{defn}

\begin{remark}
For certain classes of 2-CY-tilted algebras, such as cluster-tilted
algebras of Dynkin type, there exist algorithms to determine the
relations from their quivers~\cite{BMR06}. Moreover, the nature of
these relations allows to compute bases of paths for these algebras and
hence, thanks to Proposition~\ref{p:Tcomb}, to efficiently decide for
any vertex whether the BB-tilting module is defined.

The author has implemented these algorithms in a computer program that
determines whether two (neighboring) cluster-tilted algebras of Dynkin
type (given by their quivers) are related by a good mutation. In
particular, this allowed to verify the results of~\cite{BHL09} on a
computer.
\end{remark}

Inspired by the connectivity of the mutation graph of cluster-tilting
objects in cluster categories~\cite{BMRRT06}, the following question
naturally arises.

\begin{quest}
Let $\gL$ and $\gL'$ be two derived equivalent 2-CY-tilted algebras
(for a fixed 2-CY category $\cC$). Are they connected by a sequence of
good mutations (in $\cC$)?

In other words, does there exist a sequence of 2-CY-tilted algebras
$\gL = \gL_0, \gL_1, \dots, \gL_t = \gL'$ such that for any $0 \leq i <
t$, either $\gL_{i+1}$ is obtained from $\gL_i$ by a good mutation or
$\gL_i$ is obtained from $\gL_{i+1}$ by a good mutation?
\end{quest}

The results in~\cite{Bastian09,BHL09,BuanVatne08} show that the answer
to this question is positive for the cluster-tilted algebras of types
$A$, $\tilde{A}$ and $E$. However, there are cluster-tilted algebras of
type $D$ which are derived equivalent but not connected by a sequence
of good mutations, as shown by the following example, see also the
forthcoming paper~\cite{BHL10}.

\begin{example}
The quivers of two derived equivalent cluster-tilted algebras of type
$D_6$ are shown below.
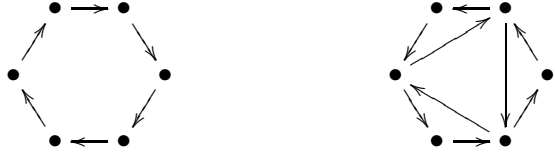
\begin{align*}
\xymatrix@=0.3pc{
& {\bullet} \ar[rr] & & {\bullet} \ar[rdd] \\ \\
{\bullet} \ar[ruu] & & & & {\bullet} \ar[ldd] \\ \\
& {\bullet} \ar[luu] & & {\bullet} \ar[ll]
}
& &
\xymatrix@=0.3pc{
& {\bullet} \ar[ldd] & & {\bullet} \ar[ll] \ar[dddd] \\ \\
{\bullet} \ar[rdd] \ar[rrruu] & & & & {\bullet} \ar[luu] \\ \\
& {\bullet} \ar[rr] & & {\bullet} \ar[ruu] \ar[llluu]
}
\end{align*}
Since these algebras are self-injective~\cite{Ringel08}, they have no
non-trivial tilting modules and in particular no good mutations.
\end{example}

\subsection{A numerical criterion for derived equivalence}
\label{ssec:numerical}

By a result of Keller and Reiten~\cite{KellerReiten07}, 2-CY-tilted
algebras are Gorenstein of dimension at most one. We now use this fact
to obtain a criterion for the derived equivalence via BB-mutation of
neighboring 2-CY-tilted algebras in terms of their Cartan matrices,
under the assumption that these matrices are invertible over $\bQ$.

For a finite dimensional algebra $A$ with $n$ non-isomorphic simples,
we denote by $I_i$ the (indecomposable) injective envelope of $S_i$.

\begin{defn}
Let $A$ be a finite dimensional algebra with $\pd_A DA < \infty$.
We define the (integral) matrix $S_A$ by the following equalities
in $K_0(\per A)$:
\begin{align*}
[I_i] = \sum_{j=1}^n \left(S_A\right)_{ij} [P_j] &,& 1 \leq i \leq n .
\end{align*}
\end{defn}

When $A$ has finite global dimension, the matrix $-S_A$ is the matrix
of the \emph{Coxeter transformation} of $A$ with respect to the basis
of $K_0(\per A)$ given by the indecomposable projectives,
see~\cite{Lenzing99}.

\begin{lemma}
Let $A$ be a finite dimensional algebra with $\pd_A DA < \infty$.
\begin{enumerate}
\renewcommand{\theenumi}{\alph{enumi}}
\item
If $\id_A A < \infty$, then $S_A$ is invertible over $\bZ$.

\item
Let $C_A$ be the Cartan matrix of $A$. Then $S_A C_A^T = C_A$.

\item
In particular, if $C_A$ is invertible over $\bQ$, then $S_A$ equals
the \emph{asymmetry} $C_A C_A^{-T}$.
\end{enumerate}
\end{lemma}
\begin{proof}
\begin{enumerate}
\renewcommand{\theenumi}{\alph{enumi}}
\item
Let $w$ be the integral matrix defined by the equalities
$[P_i] = \sum_{j=1}^n w_{ij} [I_j]$ in $K_0(\per A)$ for $1 \leq i \leq n$.
Then $w \cdot S_A$ is the identity matrix.

\item
By duality, we have $\Hom(P_l, I_i) \simeq D\Hom(P_i, P_l)$ for
$1 \leq i, l \leq n$. Thus
\begin{align*}
\sum_{j=1}^n (S_A)_{ij} (C_A)_{lj} &= \langle [P_l], [I_i] \rangle_A
= \dim_K \Hom(P_l,I_i) = \dim_K \Hom(P_i, P_l) \\
&= (C_A)_{il}
\end{align*}
and so $S_A \cdot C_A^T = C_A$.
\end{enumerate}
\end{proof}

As before, let $U = U_1 \oplus \dots \oplus U_n$ be a basic
cluster-tilting object in the 2-CY category $\cC$ and let $\gL =
\End_{\cC}(U)$ be the corresponding 2-CY-tilted algebra. Since $\gL$ is
Gorenstein, the matrix $S_\gL$ is defined. We list some of its
properties, based on the notion of \emph{index}
from~\cite{DehyKeller08}.

As the Gorenstein dimension of $\gL$ is at most
one~\cite{KellerReiten07}, there exists a minimal projective resolution
of $D\gL$ of the form
\begin{equation} \label{e:presDL}
\bigoplus_{j=1}^n P_j^{e_j} \to \bigoplus_{j=1}^n P_j^{d_j}
\end{equation}
with integers $d_j, e_j \geq 0$ for $1 \leq j \leq n$.

\begin{lemma} \label{l:djej0}
Let $1 \leq j \leq n$. Then one of $d_j$, $e_j$ vanishes.
\end{lemma}
\begin{proof}
By~\cite{KellerReiten07}, the resolution~\eqref{e:presDL} arises from a
triangle
\[
U^1 \to U^0 \to \nu U \to U^1[1]
\]
with $U^0, U^1 \in \add U$, where $\nu \simeq [2]$ is the Serre functor
in $\cC$. Since $\nu U$ is a rigid object of $\cC$, the objects $U^0$
and $U^1$ have no common summand by~\cite[Prop.~2.1]{DehyKeller08}.
This means that $d_j$ and $e_j$ cannot be both positive.
\end{proof}

It follows that one can read the terms of the minimal projective
resolution of each $I_i$ from the entries of the matrix $S_\gL$.

\begin{prop} \label{p:IpresS}
For $1 \leq i \leq n$, let
\[
\bigoplus_{j=1}^n P_j^{e_{ij}} \to \bigoplus_{j=1}^n P_j^{d_{ij}}
\]
be a minimal projective resolution of $I_i$. Then:
\begin{enumerate}
\renewcommand{\theenumi}{\alph{enumi}}
\item
For any $1 \leq j \leq n$,
\begin{align*}
d_{ij} = \begin{cases}
(S_\gL)_{ij} & \text{if $(S_\gL)_{ij} > 0$,} \\
0 & \text{otherwise;}
\end{cases}
& &
e_{ij} = \begin{cases}
-(S_\gL)_{ij} & \text{if $(S_\gL)_{ij} < 0$,} \\
0 & \text{otherwise.}
\end{cases}
\end{align*}

\item
A column of $S_\gL$ cannot contain both positive and negative entries.
\end{enumerate}
\end{prop}
\begin{proof}
Let $1 \leq j \leq n$.
We have $d_j = \sum_{i=1}^n d_{ij}$ and $e_j = \sum_{i=1}^n e_{ij}$
as sums of non-negative integers. Since one of $d_j, e_j$ vanishes
by Lemma~\ref{l:djej0}, we have that all the entries $d_{ij}$ or
all the entries $e_{ij}$ vanish.

In particular, since one of $d_{ij}, e_{ij}$ vanishes, the first
assertion follows from $(S_\gL)_{ij} = d_{ij} - e_{ij}$.

If the $j$-th column of $S_\gL$ contained entries of different signs,
this would mean that $d_{ij}$ and $e_{i'j}$ are non-zero for some $i,
i'$, which is impossible.
\end{proof}

By using the triangle $U \to U^0 \to U^1 \to U[1]$ with $U^0, U^1 \in
\add \nu U$ we obtain the dual statement.
\begin{prop}
For $1 \leq i \leq n$, let
\[
\bigoplus_{j=1}^n I_j^{d_{ij}} \to \bigoplus_{j=1}^n I_j^{e_{ij}}
\]
be a minimal injective resolution of $P_i$. Then:
\begin{enumerate}
\renewcommand{\theenumi}{\alph{enumi}}
\item
For any $1 \leq j \leq n$,
\begin{align*}
d_{ij} = \begin{cases}
(S_\gL^{-1})_{ij} & \text{if $(S_\gL^{-1})_{ij} > 0$,} \\
0 & \text{otherwise;}
\end{cases}
& & e_{ij} = \begin{cases}
-(S_\gL^{-1})_{ij} & \text{if $(S_\gL^{-1})_{ij} < 0$,} \\
0 & \text{otherwise.}
\end{cases}
\end{align*}

\item
A column of $S_\gL^{-1}$ cannot contain both positive and negative
entries.
\end{enumerate}
\end{prop}

\begin{cor} \label{c:BBcondS}
Let $k$ be a vertex without loops in the quiver of $\gL$. Then:
\begin{enumerate}
\renewcommand{\theenumi}{\alph{enumi}}
\item
$T^{\BB}_k(\gL)$ is defined if and only if $(S_\gL)_{ik} \leq 0$ for
all $1 \leq i \leq n$.

\item
$T^{\BB}_k(\gL^{op})$ is defined if and only if $(S_\gL^{-1})_{ik} \leq
0$ for all $1 \leq i \leq n$.
\end{enumerate}
\end{cor}
\begin{proof}
By Lemma~\ref{l:BBpres}, $T^{\BB}_k(\gL)$ is defined if and only if
$d_k = 0$. Now use Proposition~\ref{p:IpresS}. The proof of the second
statement is dual.
\end{proof}

\begin{thm}
Let $U$ be a cluster-tilting object in $\cC$ and let $\gL =
\End_{\cC}(U)$ and $\gL' = \End_{\cC}(\mu_k(U))$ be two neighboring
$2$-CY-tilted algebras. Assume that their quivers have no loops at $k$.

Then $\gL' = \mu^{\BB}_k(\gL)$ if and only if $(S_\gL)_{ik} \leq 0$ and
$\left(S_{\gL'}^{-1}\right)_{ik} \leq 0$ for any $1 \leq i \leq n$.
\end{thm}
\begin{proof}
Use Theorem~\ref{t:2CYtri} and Corollary~\ref{c:BBcondS}.
\end{proof}

\begin{remark}
When the Cartan matrices $C_{\gL}$ and $C_{\gL'}$ are invertible over
$\bQ$, the theorem gives an effective criterion to decide if $\gL' =
\mu^{\BB}_k(\gL)$, by examining the signs of the entries in the $k$-th
columns of the asymmetries $S_{\gL} = C_{\gL} C_{\gL}^{-T}$ and
$S_{\gL'}^{-1} = C_{\gL'}^T C_{\gL'}^{-1}$.
\end{remark}

\begin{example}
Consider the cluster-tilted algebras $\gL$ and $\gL'$ of type $D_5$
whose quivers, related by a mutation at the vertex $3$, are shown
below.
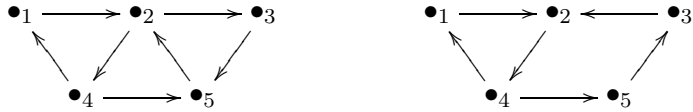
\begin{align*}
\xymatrix@=0.5pc{
{\bullet_1} \ar[rr] & & {\bullet_2} \ar[rr] \ar[ldd]
& & {\bullet_3} \ar[ldd] \\ \\
& {\bullet_4} \ar[luu] \ar[rr] & & {\bullet_5} \ar[luu]
}
& &
\xymatrix@=0.5pc{
{\bullet_1} \ar[rr] & & {\bullet_2} \ar[ldd]
& & {\bullet_3} \ar[ll] \\ \\
& {\bullet_4} \ar[luu] \ar[rr] & & {\bullet_5} \ar[ruu]
}
\end{align*}

Using the description of the relations given in~\cite{BMR06}, one
computes their Cartan matrices
\begin{align*}
C_{\gL} =
\begin{pmatrix}
1 & 0 & 0 & 1 & 0 \\
1 & 1 & 0 & 1 & 1 \\
1 & 1 & 1 & 0 & 0 \\
0 & 1 & 0 & 1 & 0 \\
0 & 1 & 1 & 1 & 1
\end{pmatrix}
& & C_{\gL'} =
\begin{pmatrix}
1 & 0 & 0 & 1 & 0 \\
1 & 1 & 1 & 1 & 1 \\
0 & 0 & 1 & 1 & 1 \\
0 & 1 & 1 & 1 & 0 \\
0 & 1 & 0 & 1 & 1
\end{pmatrix}
\end{align*}
and the corresponding asymmetries
\begin{align*}
S_{\gL} =
\begin{pmatrix}
0 & 0 & 0 & 1 & -1 \\
0 & 0 & 0 & 1 & 0 \\
1 & 0 & 0 & 0 & 0 \\
0 & 1 &-1 & 0 & 0 \\
0 & 1 & 0 & 0 & 0
\end{pmatrix}
& & S_{\gL'}^{-1} =
\begin{pmatrix}
0 & 1 &-1 & 0 & 0 \\
0 & 0 & 0 & 0 & 1 \\
0 & 0 & 0 & 1 & 0 \\
0 & 1 & 0 & 0 & 0 \\
-1 & 1 & 0 & 0 & 0
\end{pmatrix} .
\end{align*}

Since the entries in the third column of both matrices are
non-positive, the corresponding BB-tilting modules $T^{\BB}_3(\gL)$ and
$T^{\BB}_3(\gL'^{op})$ are defined, and we deduce that $\gL' \simeq
\mu^{\BB}_3(\gL)$ is obtained from $\gL$ by a good mutation at the
vertex $3$.
\end{example}

\section{Mutations of algebras of global dimension at most $2$}
\label{sec:gldim2}

We start with the following motivating example. Consider the two
skew-symmetric matrices
\begin{align*}
b = \begin{pmatrix}
0 & -1 & 0 \\
1 & 0 & -1 \\
0 & 1 & 0
\end{pmatrix}
& &
b' = \mu_2(b) = \begin{pmatrix}
0 & 1 & -1 \\
-1 & 0 & 1 \\
1 & -1 & 0
\end{pmatrix}
\end{align*}
whose associated quivers $Q$ and $Q'$ via~\eqref{e:BfromQ} are given by
those in~\eqref{e:ctaA3}. So by considering the matrices $c_Q$ and
$c_{Q'}$ defined by
\[
(c_Q)_{ij} = - \bigl| \{\text{arrows $i \to j$ in $Q$} \} \bigr| ,
\]
one can think of $b$ and $b'$ as the skew-symmetrizations $b = c_Q -
c_Q^T$ and $b' = c_{Q'} - c_{Q'}^T$.

As $Q$ is acyclic while $Q'$ is not, the path algebra $KQ$ is finite
dimensional while $KQ'$ is infinite dimensional, so they cannot be
derived equivalent. Thus, at a first stage one must introduce some
relations on $Q'$ so that by dividing its path algebra by the ideals
they generate, one obtains a finite dimensional algebra. The notion of
\emph{quivers with potentials} and their mutations, introduced
in~\cite{DWZ08}, provides a systematic way of doing this, and the
resulting algebras are known as \emph{Jacobian algebras}.

The Jacobian algebras corresponding to $Q$ and $Q'$ are precisely the
cluster-tilted algebras of Example~\ref{ex:A3}, which, despite being
closely related via near Morita equivalence, are not derived
equivalent. One approach to get a derived equivalence is to replace the
Jacobian algebras by suitable dg-algebras (the \emph{Ginzburg
algebras}). Then one is always able to interpret mutation as derived
equivalence, see~\cite{KellerYang09}.

Another approach is inspired from the fact that cluster-tilted algebras
are relation-extension of tilted algebras~\cite{ABS08}. Therefore, by
cleverly deleting arrows and interpreting them as relations in the
opposite direction, one should obtain algebras of global dimension at
most $2$, see also~\cite{ART09}, and (sometimes) interpret quiver
mutation as mutation of these algebras in the sense of
Section~\ref{sec:algmut}.

Indeed, if $A=KQ$, then by deleting the arrows in $Q'$ from $3$ to $2$
or from $2$ to $1$, one gets the algebras $\mu^-_2(A)$ and $\mu^+_2(A)$
respectively, see Example~\ref{ex:algmut}, thus interpreting the quiver
mutation at the vertex $2$ as a mutation of algebras. Moreover, one can
view $b$ and $b'$ as the skew-symmterizations of the corresponding
Euler forms,
\begin{align*}
b = c_A - c_A^T &,& b' = c_{\mu^-_2(A)} - c_{\mu^-_2(A)}^T =
c_{\mu^+_2(A)} - c_{\mu^+_2(A)}^T .
\end{align*}

\subsection{Mutations of algebras as mutations of quivers}

Let $K$ be an algebraically closed field and let $A$ be a finite
dimensional algebra of global dimension at most $2$. Denote by $Q_A$
the quiver of $A$ and recall that $c_A$ denotes the matrix of the Euler
form on $\dA$ with respect to the basis of simples. Motivated
by~\cite{ABS08,Keller09}, we make the following definition.

\begin{defn}
The \emph{extended quiver} of $A$, denoted $\wt{Q}_A$, has the same
vertices as $Q_A$, with the number of arrows from $i$ to $j$ equal to
\[
\dim_K \Ext^1_A(S_i, S_j) + \dim_K \Ext^2_A(S_j, S_i) .
\]
\end{defn}
When $A$ is hereditary, $\wt{Q}_A$ coincides with $Q_A$.

\begin{lemma}
$b_{\wt{Q}_A} = c_A - c_A^T$.
\end{lemma}
\begin{proof}
We have
\[
\begin{split}
\bigl(b_{\wt{Q}_A}\bigr)_{ij} &= -\dim_K \Ext^1_A(S_i, S_j) - \dim_K
\Ext^2_A(S_j, S_i) \\
& \quad + \dim_K \Ext^1_A(S_j, S_i) + \dim_K \Ext^2_A(S_i, S_j) \\
&= (c_A)_{ij} - (c_A)_{ji}
\end{split}
\]
since $\gldim A \leq 2$.
\end{proof}

\begin{lemma} \label{l:rQwtQ}
Let $k$ be a vertex of $\wt{Q}_A$ without loops.
\begin{enumerate}
\renewcommand{\theenumi}{\alph{enumi}}
\item
If $\mu^{-}_k(A)$ is defined and $\gldim \mu^{-}_k(A) \leq 2$, then
$r^-_k(\wt{Q}_A) = r^-_k(Q_A)$.

\item
If $\mu^{+}_k(A)$ is defined and $\gldim \mu^{+}_k(A) \leq 2$, then
$r^+_k(\wt{Q}_A) = r^+_k(Q_A)$.
\end{enumerate}
\end{lemma}
\begin{proof}
We show only the first claim, as the proof of the second is similar. By
the definition of the matrices $r^-_k$, we only need to show that
$\Ext^2_A(S_k, S_i) = 0$ for any $i \neq k$.

Let $A' = \mu^-_k(A)$ and denote by $F = \RHom(T^-_k, -)$ the
equivalence from $\dA$ to $\cD^b(A')$. From~\eqref{e:pTilt} and the
definition of $L_k$ we see that there exist $A'$-modules $M'_i$ such
that $F(S_i) \simeq M'_i$ for $i \neq k$ and $F(S_k) \simeq M'_k[-1]$
(in fact, $M'_k$ is the simple $A'$-module corresponding to $k$, but we
do not need this here). It follows that
\[
\begin{split}
\Ext^2_A(S_k, S_i) &= \Hom_{\dA}(S_k, S_i[2]) \simeq
\Hom_{\cD^b(A')}(M'_k[-1], M'_i[2]) \\
&= \Ext^3_{A'}(M'_k, M'_i) = 0
\end{split}
\]
by our assumption that $\gldim A' \leq 2$.
\end{proof}

The following proposition shows that mutations between algebras of
global dimension at most $2$ can be interpreted as mutations of the
corresponding extended quivers, as long as these quivers contain
neither loops nor $2$-cycles.

\begin{prop} \label{p:wtQmut}
Assume that $\wt{Q}_{A}$ does not contain loops nor $2$-cycles, and let
$k$ be a vertex of $\wt{Q}_A$.
\begin{enumerate}
\renewcommand{\theenumi}{\alph{enumi}}
\item \label{it:wtQminus}
If $\mu^{-}_k(A)$ is defined and $\gldim \mu^{-}_k(A) \leq 2$, then
$b_{\wt{Q}_{\mu^{-}_k(A)}} = \mu_k(b_{\wt{Q}_A})$. Moreover, if
$\wt{Q}_{\mu^{-}_k(A)}$ does not have loops or $2$-cycles, then
\[
\wt{Q}_{\mu^-_k(A)} = \mu_k(\wt{Q}_{A}) .
\]

\item \label{it:wtQplus}
If $\mu^{+}_k(A)$ is defined and $\gldim \mu^{+}_k(A) \leq 2$, then
$b_{\wt{Q}_{\mu^{+}_k(A)}} = \mu_k(b_{\wt{Q}_A})$. Moreover, if
$\wt{Q}_{\mu^{+}_k(A)}$ does not have loops or $2$-cycles, then
\[
\wt{Q}_{\mu^+_k(A)} = \mu_k(\wt{Q}_A) .
\]
\end{enumerate}
\end{prop}
\begin{proof}
By Lemma~\ref{l:rmut}, Lemma~\ref{l:rQwtQ} and
Corollary~\ref{c:rEuler},
\begin{align*}
\mu_k(b_{\wt{Q}_A}) &= \left(r^-_k(\wt{Q}_A)\right)^T b_{\wt{Q}_A}
\left(r^-_k(\wt{Q}_A)\right) \\
&= \left(r^-_k(Q_A)\right)^T b_{\wt{Q}_A} \left(r^-_k(Q_A)\right) =
\left(r^-_k(Q_A)\right)^T \left(c_A -
c_A^T\right)\left(r^-_k(Q_A)\right) \\
&= c_{\mu^-_k(A)} - c_{\mu^-_k(A)}^T = b_{\wt{Q}_{\mu^-_k(A)}} .
\end{align*}
This proves~(\ref{it:wtQminus}). The proof of~(\ref{it:wtQplus}) is
similar.
\end{proof}

\subsection{Interpretation in cluster categories}
For certain algebras $A$ with $\gldim A \leq 2$, the above results can
be refined to have an interpretation in terms of the generalized
cluster category $\cC_A$ associated with $A$ that was introduced
in~\cite{Amiot08}. Recall that $\cC_A$ is the triangulated hull of the
orbit category $\dA/\nu_A[-2]$, where $\nu_A$ denotes the Serre functor
on $\dA$.

Consider the $A$-$A$-bimodule $\Ext^2_A(DA,A)$ and the tensor algebra
$\wt{A} = T_A(\Ext^2_A(DA,A))$ known as the \emph{3-preprojective
algebra} of $A$, see~\cite{IyamaOppermann09b,Keller09}.
In this section we assume that $\wt{A}$ is finite dimensional over
$K$. It is shown in~\cite{Amiot08} that under this condition
the triangulated category $\cA$ is $\Hom$-finite and $2$-Calabi-Yau,
the image $\pi_A(A)$ of $A$ under the canonical projection
\[
\pi_A : \dA \to \dA/\nu_A[-2] \to \cA
\]
is a cluster-tilting object in $\cA$ with endomorphism algebra
\[
\End_{\cA}(\pi_A(A)) \simeq \wt{A} ,
\]
and the quiver of $\wt{A}$ is $\wt{Q}_A$.

\begin{lemma}[Amiot] \label{l:piAT}
Let $T$ be a tilting complex in $\dA$ and let $B = \End_{\dA}(T)$. If
$\gldim B \leq 2$, then $\pi_A(T)$ is a cluster-tilting object in $\cA$
whose endomorphism ring is isomorphic to $\wt{B}$.
\end{lemma}
\begin{proof}
The assumptions imply that we have a commutative diagram
\[
\xymatrix{
{\dA} \ar[rr]^{\Phi = \RHom(T,-)}_{\simeq} \ar[d]_{\pi_A}
& & {\dB} \ar[d]^{\pi_B} \\
{\cA} \ar[rr]^{\wt{\Phi}}_{\simeq} & & {\cB .}
}
\]
Since $\Phi$ maps $T$ to $B$, the isomorphism $\wt{\Phi}$ maps
$\pi_A(T)$ to $\pi_B(B)$, and the claim follows.
\end{proof}

\begin{lemma} \label{l:isoCA}
Let $X$ be a complex in $\dA$ of the form
\[
\dots \to 0 \to P^{-1} \to P^0 \to P^1 \to 0 \to \dots
\]
where $P^i$ are projectives. If $\pi_A(X) \simeq \pi_A(P)$ in $\cC_A$
for a projective $A$-module $P$, then $X \simeq P$ already in $\dA$.
\end{lemma}
\begin{proof}
For any $n \in \bZ$, denote by $\cD^{\geq n}$ the full subcategory of
$\dA$ consisting of complexes whose cohomology vanishes in degrees smaller
than $n$.
Let $F = \nu_A[-2]$. Since $\gldim A \leq 2$, we have
$F(\cD^{\geq n}) \subseteq \cD^{\geq n}$ for any $n \in \bZ$.

Since $\nu_A P$ is an injective $A$-module, we get that
$F P \in \cD^{\geq 2}$. Hence $F^m P \in \cD^{\geq 2}$, thus
$\Hom_{\dA}(X, F^m P) = 0$ for any $m > 0$.
Similarly, $F X$ is the complex of injectives
$\nu_A P^{-1} \to \nu_A P^0 \to \nu_A P^1$ concentrated in degrees $1$, $2$
and $3$, so it lies in $\cD^{\geq 1}$. It follows that
$\Hom_{\dA}(P, F^m X) = 0$ for any $m > 0$.

If $\pi_A(X) \simeq \pi_A(P)$, then by the definition of the orbit
category~\cite{Keller05}, there exist morphisms $f_m : P \to F^m X$
and $g_m : X \to F^m P$ in $\dA$ for $m \in \bZ$,
all but finitely many are zero, such that
\begin{align} \label{e:piXpiPiso}
\sum_{m \in \bZ} F^m(g_{-m}) \circ f_m = 1_P &,&
\sum_{m \in \bZ} F^m(f_{-m}) \circ g_m = 1_X .
\end{align}
But we have just shown that $f_m = 0$ and $g_m = 0$ for any $m > 0$,
so the equalities in~\eqref{e:piXpiPiso} simplify to $g_0 f_0 = 1_P$ and
$f_0 g_0 = 1_X$, giving that $X \simeq P$ in $\dA$.
\end{proof}

\begin{prop} \label{p:cAmut}
Let $k$ be a vertex of $Q_A$.
\begin{enumerate}
\renewcommand{\theenumi}{\alph{enumi}}
\item
Assume that $\mu^{-}_k(A)$ is defined and $\gldim \mu^{-}_k(A) \leq 2$.
Then $\pi_A(T^{-}_k)$ is the mutation of the canonical cluster-tilting
object $\pi_A(A)$ in $\cA$ at the vertex $k$.

\item
Assume that $\mu^{+}_k(A)$ is defined and $\gldim \mu^{+}_k(A) \leq 2$.
Then $\pi_A(T^{+}_k)$ is the mutation of the canonical cluster-tilting
object $\pi_A(A)$ in $\cA$ at the vertex $k$.
\end{enumerate}
\end{prop}
\begin{proof}
$\pi_A(T^-_k)$ is a cluster-tilting object in $\cC_A$ by
Lemma~\ref{l:piAT}. By~\eqref{e:pTilt}, it is obtained from $\pi_A(A)$
by replacing the summand $\pi_A(P_k)$ by $\pi_A(L_k)$. As $\pi_A(L_k)
\not \simeq \pi_A(P_k)$ by Lemma~\ref{l:isoCA}, we get the first claim
by the uniqueness of mutation~\cite{IyamaYoshino08}. The proof of the
second claim is similar.
\end{proof}

Thus, when $\cC_A$ is $\Hom$-finite,
by combining Proposition~\ref{p:cAmut} with~\cite[Theorem~II.1.6]{BIRSc09}
we get an alternative proof of Proposition~\ref{p:wtQmut}.

\subsection{Examples}

We give examples showing that starting with the extended quiver
$\wt{Q}_A$ of an algebra $A$ with $\gldim A \leq 2$, there are cases
where a mutation of $\wt{Q}_A$ at a vertex has two interpretations as
mutations of algebras, as well as other cases where is has no such
interpretation.

\begin{example}
Let $A^+$, $A$ and $A^-$ be the algebras given by the quivers with
relations
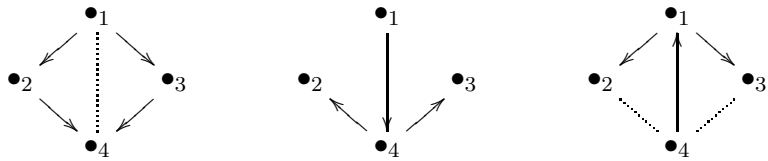
\begin{align*}
\xymatrix@=1pc{
& {\bullet_1} \ar[dl] \ar[dr] \\
{\bullet_2} \ar[dr] & & {\bullet_3} \ar[dl] \\
& {\bullet_4} \ar@{.}[uu]
} & &
\xymatrix@=1pc{
& {\bullet_1} \ar[dd] \\
{\bullet_2} & & {\bullet_3} \\
& {\bullet_4} \ar[ul] \ar[ur]
} & &
\xymatrix@=1pc{
& {\bullet_1} \ar[dl] \ar[dr] \\
{\bullet_2} & & {\bullet_3} \\
& {\bullet_4} \ar[uu] \ar@{.}[ul] \ar@{.}[ur]
}
\end{align*}
($A^+$ has a commutativity relation while $A^-$ has zero relations).
These algebras are of global dimension at most $2$ and their extended
quivers $\wt{Q}_{A^+}$, $\wt{Q}_A$ and $\wt{Q}_{A^-}$ are given by
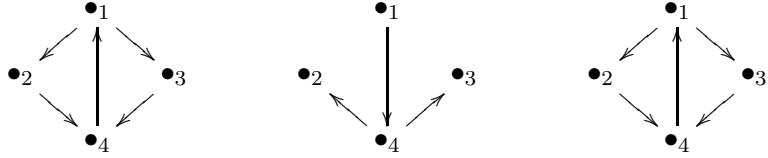
\begin{align*}
\xymatrix@=1pc{
& {\bullet_1} \ar[dl] \ar[dr] \\
{\bullet_2} \ar[dr] & & {\bullet_3} \ar[dl] \\
& {\bullet_4} \ar[uu]
} & &
\xymatrix@=1pc{
& {\bullet_1} \ar[dd] \\
{\bullet_2} & & {\bullet_3} \\
& {\bullet_4} \ar[ul] \ar[ur]
} & &
\xymatrix@=1pc{
& {\bullet_1} \ar[dl] \ar[dr] \\
{\bullet_2} \ar[dr] & & {\bullet_3} \ar[dl] \\
& {\bullet_4} \ar[uu]
}
\end{align*}

We see that $A^+ = \mu^+_4(A)$ and $A^- = \mu^-_4(A)$ and
correspondingly $\wt{Q}_{A^+} = \mu_4(\wt{Q}_A) = \wt{Q}_{A^-}$, hence
the mutation from $\wt{Q}_A$ to $\mu_4(\wt{Q}_4)$ carries two different
interpretations as mutations of algebras, leading to different, yet
derived equivalent algebras.

The passage from $A^+$ to $A^-$ can be viewed as a composition of the
following two perverse equivalences at the vertex $4$,
\[
\cD^b(A^+) \xrightarrow[\simeq]{\RHom(T^{\BB}_4(A^+),-)} \dA
\xrightarrow[\simeq]{\RHom(T^{\BB}_4(A),-)} \cD^b(A^-)
\]
arising from the BB-tilting modules $T^{\BB}_4(A^+)$ and
$T^{\BB}_4(A)$. This composition coincides with the
$2$-APR-tilt~\cite{IyamaOppermann09} of $A^+$ at the sink $4$. In fact,
one can show that any $2$-APR-tilt is a composition of the
corresponding two BB-tilts.
\end{example}

The next example shows that there
are cases where the mutation equivalence of the quivers $\wt{Q}_A$ and
$\wt{Q}_B$ can not be interpreted as a derived equivalence of the
algebras $A$ and $B$, even if the cluster categories $\cC_A$ and
$\cC_B$ are $\Hom$-finite (Claire Amiot, private communication).

\begin{example}
Let $A$ be the algebra given by the quiver with zero relations as in
the following picture.
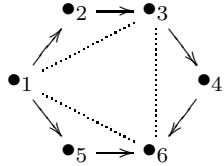
\[
\xymatrix@=0.3pc{
& {\bullet_2} \ar[rr] & & {\bullet_3} \ar[rdd] \ar@{.}[llldd] \\ \\
{\bullet_1} \ar[ruu] \ar[rdd] & & & & {\bullet_4} \ar[ldd] \\ \\
& {\bullet_5} \ar[rr] & & {\bullet_6} \ar@{.}[llluu] \ar@{.}[uuuu]
}
\]

The quiver $\wt{Q}_A$ is given by the left picture below.
If we mutate it at the vertices $2$, $4$ and $5$, we arrive at quiver $Q$
given in the right picture.
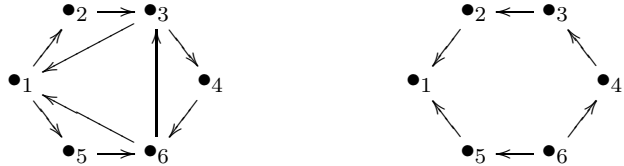
\begin{align*}
\xymatrix@=0.3pc{
& {\bullet_2} \ar[rr] & & {\bullet_3} \ar[rdd] \ar[llldd] \\ \\
{\bullet_1} \ar[ruu] \ar[rdd] & & & & {\bullet_4} \ar[ldd] \\ \\
& {\bullet_5} \ar[rr] & & {\bullet_6} \ar[llluu] \ar[uuuu]
}
& &
\xymatrix@=0.3pc{
& {\bullet_2} \ar[ldd] & & {\bullet_3} \ar[ll] \\ \\
{\bullet_1} & & & & {\bullet_4} \ar[luu] \\ \\
& {\bullet_5} \ar[luu] & & {\bullet_6} \ar[ll] \ar[ruu]
}
\end{align*}
Nevertheless, $A$ is not derived
equivalent to the algebra $B=KQ$, despite the fact that the quivers
$\wt{Q}_A$ and $\wt{Q}_B=Q$ are mutation equivalent. Indeed, even the
corresponding Euler bilinear forms $c_A$ and $c_B$ are not equivalent.

Note that the above sequence of mutations cannot be interpreted as a
sequence of perverse derived equivalences at the corresponding
vertices. Indeed, none of the complexes $T^-_2$, $T^+_2$, $T^-_4$,
$T^+_4$, $T^-_5$, $T^+_5$ is a tilting complex over $A$.
\end{example}


\def\cprime{$'$}
\providecommand{\bysame}{\leavevmode\hbox to3em{\hrulefill}\thinspace}
\providecommand{\MR}{\relax\ifhmode\unskip\space\fi MR }
\providecommand{\MRhref}[2]{%
  \href{http://www.ams.org/mathscinet-getitem?mr=#1}{#2}
} \providecommand{\href}[2]{#2}

\end{document}